\input amstex
\input xy
\input xypic
\xyoption{all}
\documentstyle{amsppt}
\document
\magnification=1200
\NoBlackBoxes
\nologo
\hoffset1.5cm
\voffset2cm
\vsize15.5cm
\def\F{\bold{F}}

\def\C{\bold{C}}

\def\Q{\bold{Q}}

\def\Z{\bold{Z}}
\def\R{\bold{R}}
\def\F{\bold{F}}
\def\N{\bold{N}}

\def\A{\Cal{A}}
\def\cM{\Cal{M}}

\def\cH{\Cal{H}}
\def\cC{\Cal{C}}
\def\cO{\Cal{O}}
\def\cV{\Cal{V}}
\def\cS{\Cal{S}}
\def\cG{\Cal{G}}
\def\cQ{\Cal{Q}}
\def\cL{\Cal{L}}
\def\cT{\Cal{T}} 
\def\cP{\Cal{P}} 
\def\cD{\Cal{D}} 
 
\def\cF{\Cal{F}}
\def\cB{\Cal{B}} 

\def\cN{\Cal{N}}

\def\cI{\Cal{I}}
\def\cJ{\Cal{J}}
\def\fm{\bold{m}}
\def\cO{\Cal{O}}
\def\cV{\Cal{V}}
\def\cB{\Cal{B}}
\def\cR{\Cal{R}}

\def\Tr{\roman{Tr}}

\def\fm{\bold{m}} 

\bigskip

\quad\quad\quad\quad  \quad\quad\quad\quad \quad\quad\quad\quad  \quad\quad\quad\quad  \quad\quad\quad\quad \quad\quad\quad\quad  {Dedicated to C. N. Yang}

\bigskip

\centerline{\bf QUANTUM OPERADS}

\medskip

\centerline{\bf No\'emie Combe,\quad  Yuri I.~Manin, \quad  Matilde Marcolli}

\bigskip

{\it ABSTRACT.}  The most standard description of symmetries of a mathematical
structure produces {\it a group}. However, when the definition of this structure is motivated  
by physics, or information theory, etc., the respective symmetry objects
might become more sophisticated: quasigroups, loops,  quantum groups, ...

In this paper, we introduce and study quantum symmetries of 
very general categorical structures: {\it operads}.

Its initial motivation were spaces of probability distributions on finite sets. 

We also investigate here how structures of quantum information, such
as quantum states and some constructions of quantum codes 
are algebras over operads. 

\medskip


\bigskip

\quad\quad\quad\quad\quad\quad\quad\quad\quad\quad\quad\quad\quad\quad {\it ... esta selva selvaggia e aspra e forte}

\quad\quad\quad\quad\quad\quad\quad\quad\quad\quad\quad\quad\quad\quad {\it che nel pensier rinova la paura ...}

\smallskip

\quad\quad\quad\quad\quad\quad\quad\quad\quad\quad\quad\quad\quad\quad\quad {\it Dante Alighieri, Inferno, Canto 1}

\bigskip

\centerline{\bf INTRODUCTION AND BRIEF SURVEY}

\medskip

 The common definition  of symmetries
of a structure given on a set $S$ (in the sense of Bourbaki) is  {\it the group of bijective maps
$S\to S$ compatible with this structure}.

\smallskip

But in fact, symmetries of various structures related to storing and transmitting inormation
such as {\it information spaces} are naturally embodied in various classes of
{\it loops} such as {\it Moufang loops}, non--associative analogs of groups
(cf. [CoMaMar21].

\smallskip

 Here are some representative examples: quasigroups, loops, and Moufang loops.
 
 \smallskip
  A {\it quasigroup}
is a set $\Cal{L}$ together with binary composition law 
$$
*: \Cal{L} \times \Cal{L} \to \Cal{L}:\quad (x_1,x_2) \mapsto x_1*x_2 =: x_3         , \eqno(0.1)
$$
such that any two elements among $\{x_1,x_2,x_3\}$ uniquely determine
the third one.

\medskip

A {\it loop} is quasigroup with two--sided identity $e\in \Cal{L}$:
it means that $e*x=x*e =x$ for any $x\in \Cal{L}$.

\medskip

Finally, a {\it Moufang loop} is quasigroup whose compostion law
satisfies the ``near--associativity'' condition
$$
(x_1*x_2)*(x_3*x_4) = x_1*((x_2*x_3)*x_4).   \eqno(0.2)
$$

\smallskip

 The idea of
symmetry embodied in a group is closely related to classical physics,
in a very definite sense, going back at least to Archimedes.
When quantum physics started to replace classical, it turned out that classical symmetries
must also be replaced by their quantum versions. For a short history
of this evolution, see pp. 1--4 of [Ma88]. As a result, the mathematical theory of {\it quantum
groups} emerged.
 
 \smallskip
 
 In this paper we suggest to apply the formalism of quantisation
 on the operadic level ([BoMa08]) to symmetries of information spaces.
 The motivation of our use of adjective ``quantum''  in [BoMa08] was sometimes too
 intuitive, but the tools developed in [Sm16], furnish a very precise
 and well axiomatized framework for this.
 
 \smallskip

  The general conception of ``Quantum Operad''
 introduced and studied here was also inspired by the introduction of
 quantum error--correcting codes with a Moufang loop action: see [CoMaMar21], Sec. 6.
 
 \smallskip
 
 We are pleased to dedicate our paper to C.  N.  Yang, who led the breakthrough
 studies of gauge symmetries in quantum field theory.
 
 \bigskip
 
 \centerline{\bf 1. QUANTUM STRUCTURES} 
 
 \smallskip
 
\centerline{\bf IN SYMMETRIC MONOIDAL CATEGORIES}
 
 \medskip
 
 {\bf 1.1. Monoidal  (= tensor) categories $\bold{V}$.} ([Sm16], Sec. 2.2, 2.3).   {\it Data:} 
 multiplication $\otimes$ of objects, with identity object $\bold{1}$ and natural isomorphisms
 $$
\alpha_{A,B;C}: (A\otimes B)\otimes C  \to A\otimes (B\otimes C), \quad
\rho_A: A\otimes \bold{1} \to A, \quad
\lambda_A: \bold{1}\otimes A \to A .
$$

{\bf 1.2. Symmetric monoidal categories.} Additional {\it twist} isomorphisms
$\tau_{A,B}: A\otimes B \to B\otimes A$, with $\tau_{A,B}\tau_{B,A} = \roman{id}_{A\otimes B}$,
plus many commutative diagrams.

\medskip

{\bf 1.3. Magmas, comagmas, bimagmas, associativity and commutativity for (co, bi)magmas
in symmetric monoidal categories.} ([Sm16], Sec. 2.4).

\smallskip

Basic data for a magma: an object $A$ with multiplication morphism $\nabla : A\otimes A \to A$.

\smallskip

Basic data fo a comagma: an object $A$ with comultiplication morhism $\Delta : A\to A\otimes A$.

\smallskip 

Basic data for a bimagma: a triple $(A,\nabla , \Delta)$ as above such that the ``bimagma diagram''
(2.4) ([Sm16], p. 49) commutes.

\smallskip

(Co, bi)--magmas in a symmetric monoidal category $\bold{V}$ are themselves objects
of respective categories. Morphisms in them are those morphisms in $\bold{V}$,
which are compatible with respective basic data.

\smallskip

{\it Unitality and counitality} structures for a magma $(A,\nabla )$   (resp. comagma $(A,\Delta)$)  are respectively the morphisms
$\eta : \bold{1}\to A$ or $\varepsilon : A \to \bold{1}$ subject to additional restrictions.

\smallskip

{\bf 1.4. Monoids, comonoids, bimonoids, and Hopf algebras in symmetric mono\-idal categories.}  ([Sm16], Def. 2.7).
They are essentially (co, bi)magmas with additional (co,bi)associativity restrictions.

\medskip

{\bf 1.5. Quantum quasigroups.} ([Sm16], Sec. 3.1). {\it A quantum quasigroup} $(A, \nabla , \Delta)$
is  bimagma, for which both {\it left composite} and {\it right composite} morphisms
are invertible:
$$
\xymatrix{A\otimes A \ar[r]^{\Delta\otimes \roman{id}_A} & \ A\otimes A \otimes A \ 
\ar[r]^{\quad \roman{id}_A \otimes \nabla} & A\otimes A},
$$
$$
\xymatrix{A\otimes A  \ar[r]^{ \roman{id}_A\otimes \Delta \quad}  & \ A\otimes A \otimes A \
\ar[r]^{\quad \nabla \otimes \roman{id}_A } & A\otimes A} .
$$
\smallskip

These morphisms are sometimes called {\it fusion operators} or {\it Galois operators}.

\smallskip

{\bf 1.6. Quantum loops.} {\it A quantum loop} in $\bold{V}$ is a biunital bigmagma $(A,\nabla ,\Delta, \eta, \varepsilon )$
such that $(A, \nabla ,\Delta )$ is a quantum quasigroup.

\smallskip

{\bf 1.7. Functoriality.} ([Sm16], Prop. 3.4). Any symmetric monoidal functor 
$$
F: (\bold{V}, \otimes , \bold{1}_{\bold{V}} ) \to ( \bold{W}, \otimes , \bold{1}_{\bold{W}})
$$
sends quantum quasigroups (resp. quantum loops) in $\bold{V}$ to quantum quasigroups
(resp. quantum loops) in $\bold{W}$.

\smallskip

{\bf 1.8. Magmas etc. in   the categories of sets with direct product.} According
to [Sm16], beginning of Sec. 3.3, in such categories comultiplication in a counital
comagma is always the respective diagonal embedding.
As a corollary, we see that quantum loops and counital quantum quasigroups in
such caregories are cocommutative and coassociative.

\smallskip

As a result, we see, that in such a category counital quantum quasigroups are equivalent
 to classical quasigroups, and quantum loops are equivalent to classical loops
 ([Sm16], Prop. 3.11).
 
 
 
\bigskip

\centerline{\bf 2. MONOIDAL CATEGORIES OF OPERADS} 

\medskip

{\bf 2.1. Graphs and their categories.} Our basic definition of graphs as quadruples 
 $(F,V,\partial, j)$ and their categories is explained in [BoMa08], Sec. 1.1, p. 251.
 There $F$, resp. $V$, are called the sets of {\it flags}, resp. {\it vertices},
 and structure maps $\partial$, resp. $j$ are called {\it boundary maps}, resp. {\it involutions}.
 Usually one flag is a pair consisting of flag as such, and a {\it label},
 that should be defined separately. {\it Geometric realization} of a graph
 is the quotient set of the disjoint union of semi--intervals $(0, 1/2]$ labeled with flags
 of this graph, modulo equivalence relation, in which $0$--points of a flag
 is glued to $1/2$ of another flag, if these flags are related by the boundary
 relation, or structure involution.
 
 \smallskip
 
 Depending on the context and/or type of labelling of $\tau$, elements of $F_{\tau}$
 might be called {\it flags, leaves, tails} ... In the study of magmatic operad
 ([ChCorGi19]) and the relevant binary trees, vertices of the relevant corollas
 are called {\it nodes}, non--root flags are called {\it left child, right child} etc.
 We will try to attach all such ``heteronyms'' to our basic terminology of 
  [BoMa08].
 
 \smallskip
Below the most typical labeling of our graphs will be (see details in [BoMa08],
Sec. 1.3.2 a) and 1.3.2 e), pp. 257--259):

\smallskip

{\it (i) Orientation.}

\smallskip

{\it (ii) Cyclic labeling.}

\smallskip

To give an orientation and cyclic labeling of corolla is essentially the same as
to define it as a {\it planar} graph: corolla, embedded into an oriented real
affine plane, with labeling compatible with its orientation.

\smallskip

Graphs endowed with various labelings form categories, upon which
the operation of disjoint union $\sqcup$ defines a monoidal structure: see [BoMa08],
Sec. 1.2.4, pp. 254--255. Our central objects of study are initially defined
only for connected graphs. Therefore, introducing this monoidal product,
we must first take care of ``empty'' (or partially empty) graphs and
explain details of their functoriality. The paper [BoMa08] is interspersed
with subsections directly or indirectly motivated by this necessity.

\smallskip

For the purposes of this paper, the most important graphs are labelled {\it trees} and
{\it forests} -- disjoint unions of trees, forming {\it ``selva selvaggia e aspra e forte''}.

\medskip

{\bf 2.2. Operads and categories of operads.}  (See [BoMa08], Sec. 1.6, p. 262). We recall here the first
definition of operads in [BoMa08], 1.6 (I), and morphisms of operads as in [BoMa08], Sec. 1.6.1.

\smallskip

First of all, we fix a symmetric monoidal category of labelled graphs $\Gamma$
with disjoint union as the monoidal structure, and a symmetric monoidal
{\it ground category} $(\cG, \otimes)$, satisfying a part of conditions 1.4 a) -- f) in [BoMa08], p. 259.

\smallskip

(i) {\it An operad is  a tensor functor between two monoidal categories
$A: (\Gamma , \sqcup )\to (\cG,\otimes )$ that sends any grafting morphism to
an isomorphism.}

\smallskip

(ii)  {\it A morphism between two operads is a functor morphism.}  

\smallskip

Denote this category of operads by $\Gamma \cG OPER$.

\medskip

{\bf 2.3.  Operads and collections as  symmetric monoidal categories.} Following [BoMa08], Sec. 1.8, we
will introduce now the monoidal ``white product'' of two operads 
$A,B : (\Gamma ,\sqcup ) \to (\cG, \otimes)$ by
the formula
$$
A\circ B (\sigma ) := A(\sigma) \otimes B(\sigma)
$$
extended to morphisms in a straightforward way.

\smallskip

Clearly, $(\Gamma \cG OPER, \circ)$ is a symmetric monoidal category.

\smallskip

An important related notion is that of {\it collection}. Starting with $\Gamma$ as above,
denote by $\Gamma COR$ its subcategory, whose objects are corollas in $\Gamma$,
and morphisms between them are isomorphisms.

\smallskip

Combining  it with the ground category $(\cG, \otimes)$ as obove, we can introduce
the category $\Gamma\cG COLL$ of $\Gamma G$--collections: its
objects are functors $A_1 :  \Gamma COR \to \cG$, and morphisms are 
natural transformations between these functors.

\smallskip

The restriction of white product $\circ$ to $\Gamma\cG COLL$ defines on it the structure
of symmetric monoidal category. If $(\cG, \otimes )$ has an identity object $\bold{1}$,
then the collection $\bold{1}_{coll}$ sending each corolla to $\bold{1}$ and each
isomorphism of corollas to the identical isomorphism of $\bold{1}$, is the 
identity collection.

\medskip

{\bf 2.4. Operads as monoids.} We briefly describe here a construction by
B. Vallette ([Val04]), reproduced in [BoMa08], Appendix A, Subsection 5.

\smallskip

We will have to use here a stronger labeling of graphs in $\Gamma$
than just orientation. Besides orientation, connected objects of $\Gamma$ 
must admit a continuous real--valued function such that it decreases
whenever one moves in the direction of orientation along each flag.
Such graphs are called {\it directed ones}
(see [BoMa08], Sec. 1.3.2 b). 
\smallskip

A graph $\tau$ is called {\it two--level graph}, if it is oriented, and if there
exists a partition of its vertices $V_{\tau} = V^1_{\tau} \sqcup V^2_{\tau}$
with the following properties:

\smallskip

(i) {\it Tails at $V^1_{\tau}$ are all inputs of $\tau$, and tails
at $V^2_{\tau}$ are all outputs of $\tau$.}

(ii) {\it All  edges in  $E_{\tau}$ go from $V^1_{\tau}$ to $V^2_{\tau}$.}

\smallskip

For any two $\Gamma G$--collections $A^1, A^2$ define their product as
$$
(A^2 \boxtimes_c A^1) (\sigma ) := \roman{colim} 
(\otimes_{v\in V^1_{\tau}} A^1(\tau_v)) \otimes (\otimes_{v\in V^2_{\tau}} A^2(\tau_v)).
$$

Here colim is taken over the category of morphisms from two level graphs to $\sigma$.

\smallskip

{\bf 2.4.1. Theorem.} {\it The product $\boxtimes_c$ is a monoidal structure on collections,
and operads are monoids in the respective monoidal category.}
\smallskip

{\bf 2.4.2. Freely generated operads.} For any $\Gamma\cG$--collection $A_1$ one can define
another collection $\cF(A_1)$ together with a canonical structure of operad on it, and for any operad
$A$ each morphism of collections $A_1 \to A$ extends to a morphism of operads
$f_A: \cF(A_1) \to A$. 

\smallskip

We can imagine $\cF(A_1)$ as the operad freely generated by the collection $A_1$.

\smallskip

{\bf 2.5. Comonoids in operadic setup.} We will now introduce a category $OP$ of operads
given {\it together with their presentations} ([BoMa08], Sec. 2.4). We start with $\Gamma$ and $\cG$
as above.

\smallskip

One object of $OP$ is a family $(A, A_1, i_A)$, where $A$ is a $\Gamma\cG$--operad,
$A_1$ is a $\Gamma\cG$--collection, such that  $f_A: \cF(A_1) \to A$ is surjective.

\smallskip

Define on $OP$ a product $\odot$ by the formula
$$
(A, A_1, i_A) \odot (B, B_1, i_B) = (C, C_1, i_C) ,
$$
in which $C_1 := A_1\circ B_1$ (cf. 2.3 above), $C$:= the minimal suboperad,
containing the image $(i_A\circ i_B)(A_1\circ B_1) \subset A\circ B$,
and $i_C$ is the restriction of $I_A\circ i_B$ on $A_1\circ B_1$.

\smallskip

{\bf 2.5.1. Theorem.} (See [BoMa08], Sec. 2.4). (i) {\it The product $\odot$ defines on $OP$ a structure
of symmetric monoidal category.}

(ii) {\it The category $OP$ is endowed with the functor of inner cohomomorphisms
$$
\underline{cohom}_{OP} : OP^{op} \times OP \to OP
$$
so that we can identify, functorially with respect to  all arguments,}
$$
\roman{Hom}_{OP} (A, C\odot B) = \roman{Hom}_{OP} (\underline{cohom}_{OP} (A,B), C)
$$

(iii) {\it Therefore, one can define canonical coassociative comultiplication
morphisms}
$$
\Delta_{A,B,C} : 
\underline{cohom}_{OP} (A, C) \to \underline{cohom}_{OP}(A,B) \odot \underline{cohom}_{OP} (B,C) .
$$
\smallskip

{\bf 2.5.2. Corollary.}  {\it For any $A$, the coendomorphism operad
$$
\underline{coend}_{OP} A := \underline{cohom}_{OP} (A,A)
$$
is a comagma in the sense of 1.3 above.}

\smallskip

{\bf 2.6. The magmatic operad.}  (See [ChCorGi19]). Below we give a brief survey
of some definitions and results from [ChCorGi19], sometimes slightly changing terminology and
notation.

\smallskip

Here objects of our basic symmetric monoidal category $(\Gamma , \sqcup )$ will be
disjoint unions of oriented trees with the following additional labeling: {\it for each tree, its outcoming flags (or leaves)
are cyclically ordered}. Corollas in it are one--vertex graphs with one root and at least
two leaves. Connected objects can be obtained from a union of disjoint corollas
by grafting each root of a corolla to one of leaves of another corolla. Morphisms
are compatible with labeling.

\smallskip

An algebra over magmatic operad is a family $(\A,*)$ consisting of a set $\A$ with binary composition law
$* : \A\times\A \to \A$.

\smallskip

Thus, corollas in the magmatic category correspond to products 
$$
(x_1*((x_2)* \dots (...(x_n)))...),
$$
and generally, connected graphs in it correspond to monomials of generic arguments with all possible
arrangements of brackets.

\smallskip

{\bf 2.7. Quasigroup monomials and planar trees.} Monomials that can be obtained
by iteration of binary multiplication $*$ as in (0.1) correspond to  {\it planar} trees:
see 2.1 above for discussion of planar corollas.
Below, discussing  quasigroups in general, and Moufang loops in particular,
we will consider connected planar trees and quasigroup monomials as encoding each other
in this way.

\bigskip

\centerline{\bf 3. MOUFANG LOOPS AND OPERADS}

\medskip

{\bf 3.1. Moufang monomials and their encoding by labeled graphs.} We will start with comparing
mathematical structures of two types: {\it labeled graphs}, and {\it Moufang monomials}.

\smallskip

The words {\it loop monomials} will refer to the following class of objects.
Let $(\Cal{L} ,*)$ be a Moufang loop in the sense of  [CoMaMar21],  Definition 5.1.1.
Let $(x_1, \dots , x_n)\in \Cal{L}$. We can produce new elements of $\cL$ from this
sequence by applying to them iterated multiplication $*$.

\smallskip

The basic examples are 
$$
x_1*x_2,  \eqno(3.1)
$$
$$
(x_1*x_2)*(x_3*x_4).   \eqno(3.2)
$$

We will encode the monomial (3.1) by a {\it cyclically labeled  oriented corolla} with one vertex
and three flags, exactly one of which is the {\it output}. The bridge from (3.1) to this
corolla might be imagined as an enrichment of it by additional labeling: $x_1*x_2$ at the output,
and $x_1$, $x_2$ at two other flags, such that $(x_1, x_2, x_1*x_2)$ corresponds to the given
cyclic labeling.

\smallskip

Since (3.2) can be obtained from (3.1) by iteration and variables change, we must explain,
how such an iteration is encoded on the level of labeled graphs. The answer
is obvious: it corresponds to {\it graftings} of certain outputs to certain inputs,
so that these outputs in the enriched picture become the inputs of the respective
iteration.

\smallskip

In this way, (3.2) becomes encoded by an {\it oriented and cyclically labeled tree},
with four ordered inputs, two edges, three vertices, and one output.

\medskip

{\bf 3.2. Passage to Moufang operad: basic identity.} According to [CoMaMar21], Def. 5.1.1,
the Moufang loops are defined as structures $(\cL,*)$ satisfying the
``near--associativity'' relations
$$
(x_1*x_2)*(x_3*x_4) = x_1*((x_2*x_3)*x_4)   \eqno(3.3)
$$
The r.h.s. of (3.3) is, in turn, encoded by an oriented and cyclically labeled tree,
with four ordered inputs (the same ones as in (3.2)), two edges, and four vertices.

\medskip

{\bf 3.3. Moufang collections.} (See [BoMa08], Sec. 1.5, pp. 259--261).
Call a {\it  Moufang corolla} an oriented cyclically ordered connected graph
with one output, and morphims are {\it isomorphisms} between them. They are objects of category $\cM COR$ (particular case of categories  $\Gamma COR$ above) Clearly, $\cM COR$ is a groupoid.

\smallskip

Choose a symmetric monoidal ground category $(\cG, \otimes )$, and define respective Moufang collections.

\medskip

{\bf 3.4. Latin square designs and their encoding by graphs.} Let $\cD = (P,L)$ be a Latin square design
as in [CoMaMar21], Def. 6.8.2.1.

\smallskip

Denote by $G^0(\cD)$ the graph, defined by the following family of data (see [BoMa08], p.251):

\smallskip

{\it Vertices} $V_{G^0(\cD)}$ are lines of $\cD$:
$$
V_{G^0(\cD)} := L.
$$
\smallskip

{\it Flags} $F_{G^0(\cD)}$ are pairs $(p,l) \in P\times L$ such that $p\in l$.

\smallskip

{\it The boundary map} $\partial_{G^0(\cD)} : F_{G^0(\cD)} \to V_{G^0(\cD)}$
sends each $(p,l)$ to $l$.

\smallskip

{\it The involution} $j_{G^0(\cD)} : F_{G^0(\cD)} \to F_{G^0(\cD)}$ sends $(p,l)$ to $(p^{\prime}, l^{\prime})$,
if $p\ne p^{\prime}$ and $l=l^{\prime}$.

\medskip

{\bf 3.4.1. Simplest examples.} Using notation from [CoMaMar21], Def. 6.8.2.1, we see that
three simplest examples correspond to  cases  $N:= \roman{card}\,L =0,1,2$.

\smallskip

The case $N=0$ is degenerate: the respective designs and graphs are empty, 
and usually are included in consideration only for categorical reasons.

\smallskip

The case $N=1$ produces a {\it corolla}: the graph with one vertex and three flags, 
and boundary map sending each flag to this vertex. The involution map is identical one.

\medskip

{\bf 3.5. From loops to Latin square designs.} Consider a ML $\cL$.
Produce three labelled  copies of $\cL$: $\cL_1$, $\cL_2$, $\cL_3$, with
pairwise empty intersections. Define the design $(P,L)$ by putting
$P:= \cL_1 \sqcup \cL_2\sqcup \cL_3$ and defining a line
as such triple of points $(x_1,x_2,x_3)$, $x_i \in \cL_i$ that
$(x_1*x_2)*x_3 = 1 \in \cL$. In this last formula we implicitly forget
labels $1,2,3$ and consider Moufang multiplication in $\cL$.

\bigskip

\centerline{\bf 4. OPERADIC STRUCTURES ON QUANTUM STATES}

\bigskip

In this section we show that the operadic structures associated to classical probabilities on finite
sets, introduced in [MarThor14], extend to {\it non--unital operads on quantum states}.

\medskip

{\bf 4.1. Operads of classical and quantum probabilities.}
We first recall the main operadic structures of classical probabilities, as introduced in [MarThor14],
which we follow for the exposition in this subsection. 
It was observed in [BFL11] and [MarThor14] that classical probabilities on finite sets come
endowed with an operad structure that describes the combination of independent subsystems.

\smallskip

Namely, denote by $\cP$ the operad  in $Sets$ with objects $\cP(n)=\Delta_n$, the simplex of probabilities
on the finite set $\{ 1, \ldots, n \}$ 
$$ 
\Delta_n=\{ P=(p_i)_{i=1}^n\,|\, p_i\geq 0, \,\, \sum_{i=1}^n p_i =1 \} ,
$$
and with composition operations
$$
 \gamma: \cP(n) \times \cP(k_1) \times \cdots \times \cP(k_n) \to \cP(k_1+\cdots + k_n) 
 $$
given by the composition of probabilities of independent subsystems
$$
 \gamma(P;P_1,\ldots, P_n)=(p_r p_{r,j_r})_{r=1,\ldots,n, j_r=1,\ldots,k_r}
  $$
for $P_r=(p_{r,j_r})_{j_r=1}^{k_r}$ and $P=(p_r)_{r=1}^n$.

\medskip

{\bf 4.1.1. Averages as an algebra over the operad $\cP$.} Making explicit in this setup
the general definitions from Sec.  2.2--2.4 above, we see that
an algebra $A$ over the operad $\cP$ in a symmetric monoidal category is a family of
morphisms
$$ 
\alpha: \cP(n) \otimes A^{\otimes n} \to A ,
$$
satisfying associativity and unitality conditions, and compatibility with 
the symmetric groups actions. 

\smallskip

The set of non--negative real numbers $\R_+$ can be seen as a category with a single object
and morphisms $x\in \R_+$ and as an algebra (in the category of small categories) 
over the operad $\cP$ with the simple operations
$$ \alpha(P; x_1,\ldots, x_n) = \sum_i p_i x_i , $$
for $P=(p_i)_{i=1}^n$ and $x_i\in \R_+$.

\medskip

{\bf 4.1.2. $A_\infty$-operad and entropy.}
As was shown in [MarThor14], there is another operadic structure in the 
setting of classical probability distributions over finite sets. 

\smallskip

Let $\cT$ be the 
$A_\infty$--operad of planar rooted tress. We say that a collection of $n$--ary information measures 
$S_n$ for $n\in \N$,
satisfies the {\it coherence condition},  if for any $n>m$, the $n$--th entropy functional
$S_n$ agrees with $S_m$, when $n-m$ of the variables are vanishing.
In other words, assume that among the probabilities $(p_1, \dots , p_n)$
the only non--zero ones are $(p_{i_1} , \dots ,  p_{i_m})$ for some $i_1 < i_2 < \dots < i_m$.
Then coherence condition means that 
$$  
 S_n (p_1, \ldots, p_n) = S_m (p_{i_1}, \ldots, p_{i_m}) . 
$$
We can now   
determine on $\R_+$ the structure of algebra over the operad $\cT$ 
with the operations  
$$
 \alpha(\tau, x_1,\ldots,x_n)= \min \{ \sum_{i=1}^n p_i x_i -\frac{1}{\beta} 
S_\tau(p_1,\ldots,p_n)\,|\, P=(p_i)\in\Delta_n \} ,
$$
where $\tau \in \cT(n)$ is a planar rooted tree with $n$ leaves, $\beta >0$ is a thermodynamic
parameter (inverse temperature). The $n$--ary entropy functional $S_\tau(p_1,\ldots, p_n)$ 
associated to the tree $\tau$ is uniquely determined by the branching structure of the tree $\tau$
and the collection of coherent entropies $S_n$. 

\medskip
{\bf 4.2. Classical probabilities from quantum states.} 
Let $\cM^{(n)}$ denote the convex set of $n\times n$--density matrices (quantum states), 
$$
 \cM^{(n)}=\{ \rho\in M_{n\times n}(\C)\,|\, \rho^*=\rho, \,\, \rho\geq 0, \,\, \Tr(\rho)=1 \} . 
 $$
 Positivity condition means here, that $\rho=a^* a$ for some $a\in  M_{n\times n}(\C)$, hence
$Spec(\rho)\subset \R_+$. In a fixed basis, the diagonal density matrices form a copy of the 
simplex $\Delta_n$ embedded in $\cM^{(n)}$.

\smallskip

There are two classical probability distributions naturally associated to a quantum state as follows.

\medskip

{\bf 4.2.1. Definition.} {\it 
Given $\rho\in \cM^{(n)}$, let 
$$
 \Lambda=(\lambda_i)_{i=1}^n \ \ \  \text{ with } \ \  \lambda_i \in Spec(\rho), 
$$
be the set of eigenvalues of $\rho$, sorted in non--increasing order, and let 
$$
 P=(p_i)_{i=1}^N \ \ \ \text{ with } \ \  p_i=\rho_{ii} 
$$
be the list of the diagonal entries of $\rho$.}

\medskip 

{\bf 4.2.2. Definition.} {\it
Given two non--increasing sequences $A=\{ a_1, \ldots, a_n \}$ and $C=\{ c_1, \ldots, c_n \}$
with $\sum_{i=1}^N a_i = \sum_{i=1}^N c_i$, one says that  $A$ majorizes $C$, or 
$A \succ C$, if for all $1\leq k\leq N$
one has $\sum_{i=1}^k a_i \geq \sum_{i=1}^k c_i$.
}

\medskip

{\bf 4.2.3. Lemma.} {\it The Shannon information of the two classical probability distributions
$\Lambda=\Lambda(\rho)$ and $P=P(\rho)$ is related by $S(P)\geq S(\Lambda)$. }

\medskip

{\it Proof.}
By Schur lemma, the sequence $\Lambda$
of eigenvalues of the hermitian matrix $\rho$ majorises the sequence $P$ of its diagonal entries, 
when both are sorted in non-increasing order.
It is well known that for probabilities $\Lambda \succ P$ is equivalent to the existence a bistochastic matrix $B$ such 
that $P =B \Lambda$. The Shannon entropy is monotonically non--decreasing under bistochastic matrices, 
so $S(B\Lambda)\geq S(\Lambda)$.  $\blacksquare$

\medskip

The eigenvalues probability $\Lambda=\Lambda(\rho)$ determines the information content
of the quantum probability $\rho$, since in the von Neumann entropy
$$ S(\rho)=\Tr(\rho \log\rho) $$
the term $\log\rho$ is defined via the spectral theorem, so that we have
$$ S(\rho)=S(\Lambda)=-\sum_i \lambda_i \log \lambda_i, $$
the Shannon entropy of the classical probability $\Lambda$.

\smallskip

We will show in the next subsections that these two classical probabilities $P(\rho)$ and
$\Lambda(\rho)$ determine two non--unital operad structures on the space of 
quantum states. The operad obtained using $P(\rho)$ has better properties and
directly agrees with the operad of classical probabilities recalled in Section~4.1 above 
when restricted to $\Delta_n \subset \cM^{(n)}$.

\medskip
{\bf 4.3. Non-unital operads.}
In the unital case, one
can equivalently describe an operad $\cO$ through the composition laws
$$
 \gamma: \cO(n) \otimes \cO(k_1)\otimes \cdots \otimes \cO(k_n) \to \cO(k_1+\cdots + k_n), 
$$
with the associativity conditions (and the symmetricity and unitality conditions in the respective cases),
or else one can describe $\cO$ through insertion operations 
$$
 \circ_i: \cO(n)\otimes \cO(m) \to \cO(n+m-1). 
$$
For $1\leq j \leq a$ and $b,c\geq 0$, with 
$X\in \cO(a)$, $Y\in \cO(b)$, and $Z\in \cO(c)$, these insertions are subject to the conditions
$$
 (X\circ_j Y)\circ_i Z = \left\{ \matrix (X\circ_i Z)\circ_{j+c-1} Y & 1\leq i < j \\
X\circ_j (Y\circ_{i-j+1} Z)& j\leq i < b+j \\
(X\circ_{i-b+1} Z)\circ_j Y & j+b \leq i \leq a+b-1.  \endmatrix \right. 
$$
The composition laws $\gamma$ satisfying the associativity condition can be obtained
from the insertions $\circ_i$ through 
$$
 \gamma(X, Y_1,\ldots, Y_n)=(\cdots (X\circ_n Y_n)\circ_{n-1} Y_{n-1}) \cdots \circ_1 Y_1). 
$$

\smallskip

While these two descriptions of operads are equivalent in the unital case, they give rise to  two
different versions of the notion of non--unital operad, see the discussion in [Markl08]. Indeed,
non--unital operads defined through the operations $\circ_i$ are also non--unital operads with the
composition operations $\gamma$, but the converse no longer holds, so the first class of 
non--unital operads is more restrictive. 

\smallskip

We will show below that the non--unital operads of quantum states belong to
the more restrictive class, according to the stronger notion of non--unital operad as in [Markl08]. 

\medskip
{\bf 4.4. The $\cQ_P$-operad of quantum states.} 
We show here that the operad $\cP$ of classical probabilities on finite sets extends
to a compatible but non--unital operad $\cQ_P$ on quantum states.

\medskip

{\bf 4.4.1. Definition.} {\it For $n\geq 1$ denote by $\cQ_P(n)=\cM^{(n)}$ the convex set of density matrices,
endowed with the composition laws
$$
 \gamma_P: \cQ_P(n)\times \cQ_P(k_1)\times \cdots \times \cQ_P(k_n)\to \cQ_P(k_1+\cdots + k_n) :
$$
$$
 \gamma_P(\rho; \rho_1,\ldots,\rho_n)= \gamma(P(\rho); \rho_1,\ldots,\rho_n) =  \pmatrix
p_1 \rho_1 & & & \\
& p_2 \rho_2 & & \\
\vdots & & \cdots & \vdots \\
& & & p_n \rho_n 
\endpmatrix  
.
$$
}

\medskip

{\bf 4.4.2. Lemma.} {\it The action of the symmetric group $\Sigma_n$ on $\cM^{(n)}$
given by $\sigma(\rho)=\sigma \rho \sigma^*$ is compatible with the action by permutation
of the coordinates on classical probabilities. It acts on the
two probability distributions $P(\rho)$ and $\Lambda(\rho)$ by
$$
 P(\sigma\rho\sigma^*)=\sigma^* P(\rho) \ \ \ \ \text{ and } \ \ \ \  \Lambda(\sigma\rho\sigma^*)=\Lambda(\rho). 
$$}

\medskip

{\it Proof.} By realising the set of classical probability distributions  $\Delta_n\subset \cM^{(n)}$
as set of diagonal density matrices in
a chosen basis, we see that $\sigma \rho \sigma^*$ permutes the entries by $\sigma^*$. The
diagonal entries of $\rho$ can be obtained as $\rho_{ii}=\Tr(\pi_i \rho)$ with $\pi_i$ the $i$-th
$1$-dimensional projection in the chosen basis, and $\Tr(\pi_i \sigma \rho \sigma^*)=\Tr(\sigma^*\pi_i \sigma\rho)
=\Tr(\pi_{\sigma^{-1}(i)} \rho)=\rho_{\sigma^{-1}(i)\sigma^{-1}(i)}$. So we have 
$P(\sigma\rho\sigma^*)=\sigma^* P(\rho)$. In the case of $\Lambda(\rho)$, since this distribution is defined
after choosing an order in which to list the eigenvalues of $\rho$, such as non--increasing order,
we have $\Lambda(\rho)=\Lambda(\sigma\rho\sigma^*)$, since both matrices have the same spectrum.
$\blacksquare$

\medskip

{\bf 4.4.3. Proposition.} {\it The convex sets $\cQ_P(n)$ with the composition operations $\gamma_P$ of 
Definition~4.4.1 determine a non-unital symmetric operad $\cQ_P$ that restricts to the unital operad
$\cP$ on classical probabilities $\Delta_n \subset \cM^{(n)}$. }

\medskip

{\it Proof.}
It is clear by the definition of the composition operations $\gamma$ that they agree with the
composition operations of the operad $\cP$ when restricted to classical probabilities $\Delta_n\subset \cM^{(n)}$.
We need to check that they satisfy the associativity and symmetry axioms on the larger set $\cM^{(n)}$ of
quantum states.

The associativity condition is given by the identities
$$ 
\gamma( \gamma(\rho^{(m)}; \rho^{(n_1)}, \ldots, \rho^{(n_m)}); \rho^{(r_{1,1})}, \ldots, \rho^{(r_{1,n_1})}, \ldots, \rho^{(r_{m,1})}, \ldots, \rho^{(r_{m,n_m})}) = 
$$
$$
 \gamma(\rho^{(m)}; \gamma(\rho^{(n_1)}; \rho^{(r_{1,1})}, \ldots, \rho^{(r_{1,n_1})}), \ldots, \gamma(\rho^{(n_m)}; \rho^{(r_{m,1})}, \ldots, \rho^{(r_{m,n_m})})), 
 $$
for $\rho^{(m)}\in \cQ(m)$, $\rho^{(n_i)}\in \cQ(n_i)$, $i=1,\ldots, m$, 
and $\rho^{(r_{i,\ell_i})} \in \cQ(r_{i,\ell_i})$ with $\ell_i=1,\ldots, n_i$. The left--hand--side is
$$
 \gamma \left( \pmatrix
\rho^{(m)}_{11} \rho^{n_1} & & & \\
& \rho^{(m)}_{22} \rho^{n_2} & & \\
\vdots & & \cdots & \vdots \\
& & & \rho^{(m)}_{mm} \rho^{n_m}
\endpmatrix ; \rho^{(r_{1,1})}, \ldots, \rho^{(r_{1,n_1})}, \ldots, \rho^{(r_{m,1})}, \ldots, \rho^{(r_{m,n_m})}\right) = 
$$
$$
   \pmatrix 
\rho^{(m)}_{11} \rho^{n_1}_{11}  \rho^{(r_{1,1})} &  & & & & \\
& \ddots & & & & & \\
& &  \rho^{(m)}_{11} \rho^{n_1}_{n_1 n_1}  \rho^{(r_{1,n_1})}  & & & & \\
& & &  \ddots & &    &  \\
& & & & \rho^{(m)}_{mm} \rho^{n_1}_{11}  \rho^{(r_{m,1})} & & \\
& & & & & \ddots & \\
& & & & & &\rho^{(m)}_{mm} \rho^{n_1}_{n_m n_m}  \rho^{(r_{m,n_m})} 
\endpmatrix 
,
$$
which agrees with the right--hand--side
$$
 \gamma\left( \rho^{(m)}; \pmatrix \rho^{n_1}_{11}  \rho^{(r_{1,1})} & & \\ & \ddots & \\ & & \rho^{n_1}_{n_1 n_1}  \rho^{(r_{1,n_1})} \endpmatrix, \ldots , 
\pmatrix \rho^{n_m}_{11}  \rho^{ ( r_{m,1} ) } & & \\ & \ddots & \\ & & \rho^{n_m}_{n_m n_m}  \rho^{ ( r_{m,n_m} ) } \endpmatrix \right) .
$$
The compatibility with the symmetric group action for $\cQ_P$ is obtained directly from Lemma~4.4.2. Indeed,
symmetric property of an operad is expressed by 
the following two identities, for permutations $\sigma_i\in \Sigma_{n_i}$ and $\sigma\in \Sigma_m$.
The first condition is
$$
 \gamma_P(\sigma(\rho); \rho_{\sigma^{-1}(1)},\ldots,\rho_{\sigma^{-1}(m)}) =
\tilde\sigma(\gamma_P(\rho; \rho_1,\ldots,\rho_m)), 
$$
where on the right-hand-side $\tilde\sigma \in \Sigma_{n_1+\cdots + n_m}$ is the
permutation that splits the set of indices into blocks of $n_i$ indices and permutes the
blocks by $\sigma$. The second symmetric group condition is
$$
 \gamma_P(\rho; \sigma_1(\rho_1),\ldots,\sigma_m(\rho_m)) = \hat\sigma(\gamma_P (\rho; \rho_1,\ldots,\rho_m)), 
$$
where on the right--hand--side $\hat\sigma  \in \Sigma_{n_1+\cdots + n_m}$ is the permutation that
acts on the $i$--th block of $n_i$ indices as the permutation $\sigma_i$.

In the first case, we have
$$
 \gamma_P(\sigma(\rho); \rho_{\sigma^{-1}(1)},\ldots,\rho_{\sigma^{-1}(m)})=\gamma_P (\sigma^{-1} P(\rho);
\rho_{\sigma^{-1}(1)},\ldots,\rho_{\sigma^{-1}(m)}) , 
$$
which is the same as $\tilde\sigma \gamma_P(\rho; \rho_1,\ldots,\rho_m)) \tilde\sigma^*$ in $\cM^{(n_1+\cdots+n_m)}$.
A similar argument proves the second relation.

The operad is non--unital. Indeed, the unit axiom is only satisfied for
$\rho=1\in \cQ(1)$ with $\gamma_P(1;\rho)=\rho$,
but it fails when $\rho_i= 1\in \cQ(1)$, where the composition gives instead 
$\gamma_P(\rho;1,\ldots,1)=P(\rho)$. 
The unit axiom $\gamma_P(\rho;1,\ldots,1)=\rho$ 
is satisfied on the subset $\Delta_n \subset \cM^{(n)}$ of classical probabilities,
where the operad agrees with the unital operad $\cP$.
$\blacksquare$

\medskip

{\bf 4.4.4. Proposition.} {\it The composition laws $\gamma_P$ of the non-unital operad $\cQ_P$
are induced by insertion operations $\circ_i: \cQ_P(n)\times \cQ_P(m) \to \cQ_P(n+m-1)$, hence
the operad $\cQ_P$ is also a non--unital operad in the stronger sense.}

\medskip

{\it Proof}. 
For the composition operations $\gamma_P$ to be obtained from insertions $\circ_i$,
we need to have 
$$
 \gamma_P(\rho; \rho_1,\ldots,\rho_n) =(\cdots (\rho\circ_n \rho_n)\cdots \circ_1 \rho_1), 
$$
for $\rho\in \cM^{(n)}$ and $\rho_i\in \cM^{(n_i)}$. We can define morphisms  
$$
 \circ_i : \cQ_P(n) \times \cQ_P(m) \to \cQ_P(m+n-1) 
 $$
as the operations $(\rho,\rho')\mapsto \rho\circ_i \rho'$ that take a density matrix $\rho\in \cM^{(n)}$
and replace the $i$--th row and columns with $m$ rows and $m$ columns, respectively, where
all the entries outside of the $m\times m$--block around the diagonal are zero, and the
diagonal block is given by the matrix $\rho_{ii} \rho'$. Clearly this implies that the
repeated application of these insertions performed in the order
$(\cdots (\rho\circ_n \rho_n)\cdots \circ_1 \rho_1)$ produces exactly the matrix 
$\gamma_P(\rho; \rho_1,\ldots,\rho_n)$.
$\blacksquare$ 

\medskip

{\bf 4.5. The $\cQ_\Lambda$--operad of quantum states.} 
We can construct, in a very similar way, another operad of quantum states, using the
classical probabilities $\Lambda(\rho)$ instead of $P(\rho)$. The resulting operad has
slightly different properties, coming from the choice of an ordering of the eigenvalues.

\medskip

{\bf 4.5.1. Definition.} {\it For $n\geq 1$ let $\cQ_\Lambda(n)=\cM^{(n)}$ the convex set of density matrices,
endowed with the composition laws
$$ \gamma_\Lambda: \cQ_\Lambda(n)\times \cQ_\lambda(k_1)\times \cdots \times 
\cQ_\Lambda(k_n)\to \cQ_\Lambda(k_1+\cdots + k_n) $$
$$ \gamma_\Lambda(\rho; \rho_1,\ldots,\rho_n)= \gamma(\Lambda(\rho); \rho_1,\ldots,\rho_n) =  \pmatrix
\lambda_1 \rho_1 & & & \\
& \lambda_2 \rho_2 & & \\
\vdots & & \cdots & \vdots \\
& & & \lambda_n \rho_n 
\endpmatrix $$
with $\lambda_i$ the eigenvalues of $\rho$ listed in non-increasing order. 
}

\medskip

{\bf 4.5.2. Proposition.} {\it The convex sets $\cQ_\Lambda(n)$ with the composition operations $\gamma_\Lambda$ of 
Definition~4.5.1 determine a non-unital non-symmetric operad $\cQ_\Lambda$. The composition laws $\gamma_\Lambda$ are induced by insertion operations $\circ_i: \cQ_\Lambda(n)\times \cQ_\Lambda(m) \to \cQ(n+m-1)$.}

\medskip

{\it Proof.}
The argument is completely analogous to Propositions~4.4.3 and 4.4.4. The associativity requirement for
the composition laws $\gamma_\Lambda$ follows as in Propositions~4.4.3, using the fact that 
$$
 Spec\, \pmatrix
\lambda_1 \rho^{n_1} & & & \\
& \lambda_2 \rho^{n_2} & & \\
\vdots & & \cdots & \vdots \\
& & & \lambda_m \rho^{n_m}
\endpmatrix = \bigcup_i \lambda_i \, Spec(\rho^{n_i}). 
$$
The operad $\cQ_\Lambda$ is non--unital for the same reason as $\cQ_P$, namely
$\gamma_\Lambda(\rho;1\ldots,1)=\Lambda(\rho)\neq \rho$. Even in the case of
diagonal matrices these differ in general by a permutation.

The operad is non--symmetric, because
we must choose an
ordering of  the eigenvalues in $\Lambda(\rho)$, for example, non--increasing ordering. 
This breaks the symmetric group action in the first of the two identities, making
$\cQ_\Lambda$ non--symmetric. 

The insertion operations $\circ_i$ are as in Proposition~4.4.4, but 
the central $m\times m$ block is now of the form $\lambda_i \rho'$
with $\lambda_i$ the $i$--th eigenvalue of $\rho$ in non-increasing order.
$\blacksquare$

\medskip
{\bf 4.6. Trees of projective quantum measurements.} We consider the range of
compositions of insertion maps of the operad $\cQ_P$ and quantum
channels associated to these ranges.

\smallskip

First observe that the  image of the insertion map
$$ \circ_i : \cQ_P(n)\times \cQ_P(m) \to \cQ_P(n+m-1) $$
consists of the set of those density matrices $\rho \in \cM^{(n+m-1)}$ that are
block diagonal with one $(n-1)\times (n-1)$--block and one $m\times m$--block. 
Moreover, all block diagonal density matrices are in the image of some composition of
insertion maps. These are quantum states that decompose nontrivially into
disjoint states with orthogonal ranges. On the other hand, a block diagonal
density matrix can be obtained in more than one way through a composition of
insertion maps. 

\smallskip

Operators on density matrices are described by quantum channels,
namely quantum measurements realized by completely positive maps.
A particular class of such operators consists of {\it projective
quantum measurements}.

\smallskip

In the following we work with  a finite dimensional
Hilbert space $\cH$ of dimension $N$ with a chosen orthonormal basis,
and we denote as before by $\cM^{(N)}$ the set of density matrices, 
written in the chosen basis.

\medskip

{\bf 4.6.1. Definition.} {\it
A projective quantum measurement is a family $\Pi=\{ \Pi_i \}_{i=1}^n$ of projectors
$\Pi_i^*=\Pi_i=\Pi_i^2$ in a finite dimensional Hilbert space $\cH$ of dimension
$N$, that are mutually orthogonal, $\Pi_i \Pi_j =\delta_{ij} \Pi_i$,
and satisfy the comdition $\sum_i \Pi_i =1$. The outcome of the projective measurement $\Pi$
on a quantum state given by  density matrices $\rho\in \cM^{(N)}$ is
$$ \rho_i=\frac{\Pi_i \rho \Pi_i}{\Tr(\Pi_i \rho)} \ \ \text{ with probability } \ \ p_i=\Tr(\Pi_i \rho). $$
The projective quantum channel $\Pi$ then maps
$\rho\mapsto \Pi(\rho)=\sum_i p_i \rho_i$. }

\medskip

The range in $\cM^{(N)}$ of a composition of insertion maps is specified by assigning
a decomposition $N=k_1+\cdots+k_n$ and the locus $\cM_{k_1,\ldots,k_n}\subset \cM^{(N)}$
of density matrices that are block--diagonal in the chosen basis, with $n$ blocks of size $k_i$.
The following is immediate by construction.

\medskip

{\bf 4.6.2. Lemma.} {\it Let $\Pi$ be
the projective measurement $\Pi=\{ \Pi_i \}_{i=1}^n$, where $\Pi_i$ is the orthogonal projection
onto the span of the $i$--th subset of $k_i$ basis elements, given by a quantum
channel. It  maps $\Pi: \cM^{(N)} \to \cM_{k_1,\ldots,k_n}$, assigning to $\rho$ the block--diagonal
density matrix $\Pi(\rho)=\sum_i p_i \rho_i$. }

\medskip

Just as elements of $\cM_{k_1,\ldots,k_n}$ can be obtained in different ways as
compositions of insertion operations of the operad $\cQ_P$ of quantum states,
with different tree structures, the quantum channel $\Pi:\cM^{(N)} \to \cM_{k_1,\ldots,k_n}$
can also be realized through different tree structures.

\smallskip

Let $\tau$ be a planar rooted tree with $n$ leaves labeled by 
the non--negative integers $k_i$. We view the tree $\tau$ as oriented 
from the leaves toward the root. We assign to the root vertex $v_0$,
the identity projector $\Pi^{(v_0)}=1$. Let $v$ be any vertex
in the tree, and consider the set of incoming edges $e$ at $v$, with $v_e=s(e)$
their source vertices. Let $\{ \Pi^{(v_e)} \}_{t(e)=v}$ be a set of orthogonal
projections with $\sum_{e: t(e)=v} \Pi^{(s(e))}=\Pi^{(v)}$. Let $\Pi_i$ denote the
resulting projectors at the leaves. We write, for $t(e)=v$,
$$
 \rho^{(s(e))} =
\frac{ \Pi^{ ( s(e) ) } \rho^{ (v) } \Pi^{ ( s(e) ) } } { \Tr ( \Pi^{ (s(e)) } \rho^{ (v) } ) }, 
$$
with $\rho^{(v_0)}=\rho$. 

\smallskip

Consider the quantum measurement $\Pi_\tau$ that
assigns to a density matrix
$\rho\in \cM^{(N)}$, with $N=k_1+\cdots+k_n$ 
outcomes $\rho_i^\tau$ with probabilities $p_i^\tau$ where 
$$ p_i^\tau=\prod_w \Tr(\Pi^{(w)}_{i_w} \rho^{(w)}) $$
with the product over the vertices on the directed path from the $i$--th leaf to the root of $\tau$, with $i_w$
indicating the direction at the vertex $w$ along this path and the
$\rho_i^\tau$ are obtained by repeatedly computing $\rho^{(s(e))}$ from $\rho^{(t(e))}$ along
the path connecting the $i$--th leaf to the root. 

\medskip

{\bf 4.6.3. Lemma.} {\it 
All the quantum channels $\Pi_\tau$ obtained in this way, are just the same
quantum channel $\Pi:\cM^{(N)} \to \cM_{k_1,\ldots,k_n}$.}

\medskip

{\it Proof.}
This can be seen very easily by writing
$$ \prod_w \Tr(\Pi^{(w)}_{i_w} \rho^{(w)}) 
= \prod_\ell \frac{ \Tr(  \Pi^{ (w_\ell) } \Pi^{ (w_{\ell-1}) }
\cdots \Pi^{ (w_0) } \rho) } 
{ \Tr( \Pi^{ (w_{ \ell - 1 }) } \cdots \Pi^{ (w_0) } \rho ) } 
= \prod_\ell \frac{\Tr(\Pi^{(w_\ell)} \rho)}{\Tr(\Pi^{(w_{\ell-1})} \rho)}
= \Tr(\Pi_i \rho) $$ 
with $\Pi_i$ the projection at the leaf, since we have
$\Pi^{(s(e))}_j \Pi^{(t(e))}_e =\Pi^{(s(e))}_j$ as $\Pi^{(s(e))}_j$ is a projection 
onto a subspace of the range of $\Pi^{(t(e))}_e$. We similarly obtain 
$\rho^\tau_i=\rho_i=\frac{\Pi_i \rho \Pi_i}{\Tr(\Pi_i \rho)}$.
$\blacksquare$

\medskip

{\bf 4.7. Entropy functionals.} 
Consider a family of quantum entropy functionals
$$ S_n: \cM^{(n)} \to \R $$
satisfying the consistency condition that $S_n$ restricts to $S_k$ for $k<n$ over
any copy of $\cM^{(k)}$ embedded in $\cM^{(n)}$ as density matrices with a set of
$n-k$ vanishing eigenvalues. 

\smallskip

Examples of such consistent collections of entropy functionals include the
von Neumann entropy 
$$ \cN(\rho)=-\Tr(\rho \log \rho), $$
or, for a real parameter $q>0$ with $q\neq 1$,
the quantum R\'enyi entropy 
$$ Ry_q(\rho)=\frac{1}{1-q} \log \Tr(\rho^q) $$
and the quantum Tsallis entropy
$$ Ts_q(\rho)=\frac{1}{1-q} (\Tr(\rho^q)-1). $$

\smallskip

We obtain a family of quantum entropies associated to trees in the following way, as
a direct generalization of the entropy functionals $S_\tau$ for classical
probabilities constructed in [MarThor14].

\medskip

{\bf 4.7.1. Proposition.} {\it 
A tree $\tau$ with $n$ leaves labeled by integers $k_i\geq 1$, 
together with a coherent family $\{ S_n \}$ of
quantum entropies, determines an entropy functional
$$ 
S_\tau: \cM^{(N)} \to \R 
$$
for $N=k_1+\cdots+k_n$.}

\medskip

{\it Proof.}
In the case where $\tau$ is a corolla with a single root vertex and $n$ leaves,
we set
$$ 
S_\tau(\rho)=S(P)+\sum_i p_i S(\rho_i) 
$$
with $p_i=\Tr(\Pi_i \rho)$, resp. $\rho_i=\frac{\Pi_i \rho \Pi_i}{\Tr(\Pi_i \rho)}$, are
the probabiities, resp. density matrices, at the leaves. In the case of the
von Neumann entropy, by the extensivity property, this is the same
as $\cN(\Pi^\tau(\rho))=\cN(\sum_i p_i \rho_i)$. Inductively assuming
that $S_\tau$ is constructed for all trees with less than $n$ leaves,
consider the subtrees $\tau_j$, $j=1,\ldots, m$, attached at the root vertex $v_0$,
each with a set $L_j$ of leaves, $\roman{card}\, L_j <n$. By the construction of the
quantum channel $\Pi^\tau$ we have a system $\Pi_j$ of orthogonal projections
with $\sum_j \Pi_j=1$ associated to the incoming edges $e_j$ at the root vertex,
and probabilities $p_j=\Tr(\Pi_j \rho)$ and density matrices $\rho_j=\frac{\Pi_j \rho \Pi_j}{\Tr(\Pi_j \rho)}$
associated to the root vertices $v_j$ of the subtrees $\tau_j$. We then set
$$
 S_\tau(\rho)=S(P)+\sum_j p_j S_{\tau_j}(\rho_j), 
 $$
where  $S(P)$ is the Shannon entropy of the classical probability $P=(p_j)$ and
$S_{\tau_j}$ are the entropy functionals inductively constructed for the subtrees $\tau_j$
with less than $n$ leaves. This suffices to determine $S_\tau$ uniquely. 
$\blacksquare$

\medskip

{\bf 4.7.2. Remark.} In the case of the von Neumann entropy, the extensivity property and the
identification of all the quantum channels $\Pi^\tau$ with the 
quantum channel $\Pi:\cM^{(N)} \to \cM_{k_1,\ldots,k_n}$, seen in Lemma~4.6.3 above, 
imply that $\cN_\tau(\rho)=\cN(\sum_i p_i \rho_i)=S(P)+\sum_i p_i \cN(\rho_i)$
for all $\tau$. This is not the case for non--extensive entropies like R\'enyi and Tsallis.

\medskip

{\bf 4.8. $A_\infty$--operad of quantum channels.} 
We now investigate structures on
quantum states that generalize the operations based on classical probabilities
that we recalled in Section~4.1.2. 

\medskip

{\bf 4.8.1. Definition.} {\it 
A quantum channel $\Phi: \cM^{(N)}\to \cM^{(N)}$ is a trace preserving completely positive map.
It is well known that any such map can be represented, in a non-unique way, in Kraus form,
namely as
$$ \Phi(\rho)=\sum_i A_i \rho\, A_i^* $$
for a collection $\{ A_i \}$ of operator satisfying the condition $\sum_i A_i^* A_i=1$.}

\medskip

The projective quantum channels considered in the previous subsections are
those for which the operators $A_i$ are mutually orthogonal projectors.

\smallskip

We can construct more general quantum channels associated to rooted trees.
We consider a finite dimensional complex Hilbert space $\cH$ of dimension $N$. All
operators here will be linear operators on $\cH$.

\medskip

{\bf 4.8.2. Definition.} {\it
Let $\tau$ be a planar rooted tree with $n$ leaves. 
We consider $\tau$ oriented from the leaves to the root.
A tree quantum channel $C^\tau_A$ is an assignment of
operators $A=\{ A_e \}_{e\in E(\tau)}$ to the edges of $\tau$, satisfying  
the condition that at each vertex $v$
$$ \sum_{e\,:\, t(e)=v} A_e^* A_e = 1. $$
}

\medskip

{\bf 4.8.3. Lemma.} {\it The tree quantum channels $C^\tau_A$ of Definition~4.8.2 are
quantum channel as in Definition~4.8.1, acting 
on density matrices $\rho\in \cM^{(N)}$ by
$$ 
C^\tau_A(\rho)=\sum_{i=1}^n A_{e_{i,1}}\cdots A_{e_{i,m_i}} \rho\,  A_{e_{i,m_i}}^*\cdots A_{e_{i,1}}^*,
$$
where the sum is over the leaves of $\tau$ and for each $i=1,\ldots, n$
we consider the oriented path $e_{i,1},\ldots, e_{i,m_i}$ from the $i$--th
leaf to the root, $s(e_{i,1})=v_i$, $t(e_{i,j})=s(e_{i,j+1})$, $t(e_{i,m_i})=v_0$,
the root vertex. }

\medskip

{\it Proof.} This is essentially the Kraus form of a quantum channel, since we have
$$
 \sum_{i=1}^n A_{e_{i,m_i}}^*\cdots A_{e_{i,1}}^*\, A_{e_{i,1}}\cdots A_{e_{i,m_i}} =1. 
$$
Indeed, write $A_{v}:=A_{e_1}\cdots A_{e_m}$ for the composition of the operators $A_e$
 along the oriented path from the vertex $v$ to the root $v_0$. So we write the
 above as $\sum_i A_i^* A_i$. Starting at the leaves and considering the adjacent
 vertices, we can rewrite the sum as 
$$ 
\sum_v \sum_{i\,:\, t(e_i)=v} A_v^* A_{e_{i,1}}^*\, A_{e_{i,1}} A_v =\sum_v A_v^* A_v ,
$$
where the set $\{ i\,:\, t(e_i)=v \}$ is non--empty only for the vertices $v$ adjacent to the leaves.
This reduces by one the length of the path. Thus, we obtain inductively that the normalisation
$\sum_v A_v^* A_v=1$ holds, with the condition $\sum_{e\,:\, t(e)=v} A_e^* A_e = 1$ implying 
that it holds for length one. 
$\blacksquare$

\medskip

{\bf 4.8.4. Theorem.} {\it The tree quantum channels $C^\tau_A$ form an $A_\infty$--operad $\cQ\cC$. }

\medskip

{\it Proof.} 
Consider $\Z$--modules $\cQ\cC(n) := span_\Z \{ C^\tau_A\,|\, \tau \in \cT(n) \}$ where $\cT(n)$ is the $A_\infty$-operad
of planar rooted trees. The operadic composition laws 
$$
 \gamma_{\cQ\cC}: \cQ\cC(n) \otimes \cQ\cC(k_1)\otimes \cdots \otimes \cQ\cC(k_n) \to \cQ\cC(k_1+\cdots +k_n) 
$$
are given by
$$
 \gamma_{\cQ\cC} (C^\tau_A; C^{\tau_1}_{A_1}, \ldots, C^{\tau_n}_{A_n}) =
C^{\gamma_{\cT}(\tau;\tau_1,\ldots,\tau_n)}_{A\cup \{A_1,\ldots,A_n\}}.
$$
The associativity, unity, and symmetric properties of $\cQ\cC$ follow directly from the same properties
of the operad $\cT$. The DG structure of $\cQ\cC$ is also inherited from the DG--structure of the
$A_\infty$--operad $\cT$, with the differential given by edge contractions
$$ 
d C^\tau_A = \sum_{\tau'\,:\, \tau=\tau'/e} \epsilon\, C^{\tau'}_{A'}, 
$$
where $\epsilon =(-1)^{\ell(e)}$ with $\ell(e)$ the number of edges below and to the left of $e$ in $\tau$
with respect to the planar structure. 

\smallskip

The collection of operators $A'$ on a tree $\tau'$ with $\tau=\tau'/e$ agrees with $A$
on all edges $e'$ with $t(e')\neq t(e), s(e)$ and is of the following form on the remaining
edges. Let $E_t$ be the set of edges $e'$ of $\tau$ such that $t(e')=t(e)$ in $\tau'$ and $E_s$ the set of
edges $e'$ in $\tau$ such that $t(e')=s(e)$ in $\tau'$. The set $E_t\cup E_s$ consists of all the 
edges of $\tau$ with the same target vertex $v$ in $\tau$ that is split into two vertices $s(e), t(e)$
in $\tau'$. Thus, in $C^\tau_A$ we have the relation $\sum_{e' \in E_t\cup E_s} A_e^* A_e=1$. 

\smallskip

Let $B_t:=\sum_{e' \in E_t} A_{e'}^* A_{e'}$ and $B_s:=\sum_{e' \in E_t} A_{e'}^* A_{e'}$. These
are positive operators, namely $\langle B v,v\rangle\geq 0$ for all $v\in\cH$. Put $N_s=\roman{card}\, E_s$
and consider the operators $\frac{1}{N_s} B_t$ and $A_{e'}^* A_{e'}+\frac{1}{N_s} B_t$ for
$e' \in E_s$. These are also positive operators. Thus, we can write
$B_s=A^* A$ and $A_{e'}^* A_{e'}+\frac{1}{N_s} B_t =\tilde A_{e'}^* \tilde A_{e'}$ for
some operators $A$ and $\tilde A_{e'}$. 

\smallskip

We then take $A'_e:= A$ and $A'_{e'}:=\tilde A_{e'}$ for $e'\in E_s$. This completes the
description of $A'$ on $\tau'$ in a way that still satisfies the conditions at vertices
$$
 \sum_{e'\,:\, t(e')=s(e)} {A'_{e'}}^* A'_{e'} = \sum_{e'\in E_s} (A_{e'}^* A_{e'}+\frac{1}{N_s} B_t)=
\sum_{e'\in E_s} A_{e'}^* A_{e'} +\sum_{e'\in E_t} A_{e'}^* A_{e'} =1,
$$
$$
 \sum_{e'\,:\, t(e')=t(e)} {A'_{e'}}^* A'_{e'} ={A'_e}^* A'_e + \sum_{e'\in E_t} A_{e'}^* A_{e'} =
\sum_{e'\in E_s} A_{e'}^* A_{e'} +\sum_{e'\in E_t} A_{e'}^* A_{e'} =1 . 
$$

\smallskip

Note that, in order to make the composition operations 
compatible with the differential of the DG-structure, they should also include appropriate 
signs, as specified for instance in [Vor01]. We will omit the details as they
are exactly the same as in the original case of the $A_\infty$--operad $\cT$.
$\blacksquare$

\medskip

Taking formal linear combinations of quantum channels, as in the definition of $\cQ\cC$ above, 
has the advantage of being able to define the differential and DG-structure described in 
Theorem~4.8.4. However, it is somewhat unnatural, since the positivity property of quantum channels
is lost in the linear combinations. It is therefore more natural in this context to consider only convex 
combinations. This leads to a variant of the operad $\cQ\cC$ of tree quantum channels. 

\medskip

{\bf 4.8.5. Definition.} {\it Let $\cQ\cC^+$ be given by
$\cQ\cC^+(n)=\text{convex span}\{ C^\tau_A\,|\, \tau \in \cT(n) \}$.
The $\cQ\cC^+(n)$ are convex sets rather than $\Z$-modules (or vector spaces).
The composition operations are the same as the $\gamma_{\cQ\cC}$ (without signs),
$$ \gamma_{\cQ\cC^+} (C^\tau_A; C^{\tau_1}_{A_1}, \ldots, C^{\tau_n}_{A_n}) =
C^{\gamma_{\cT}(\tau;\tau_1,\ldots,\tau_n)}_{A\cup \{A_1,\ldots,A_n\}}. $$
}

\medskip

Then $\cQ\cC^+$ defined in this way is an operad, though it no longer has a DG-structure. 

\medskip

{\bf 4.8.6. Proposition.} {\it The convex set $\cM^{(N)}$ has the structure of an algebra 
over the operad $\cQ\cC^+$.}

\medskip

{\it Proof.} 
The operations $\alpha: \cQ\cC^+(n) \otimes {\cM^{(N)}}^{\otimes n} \to \cM^{(N)}$ are given by
a slight modification of the action of tree quantum channels of Lemma~4.8.3, 
$$
 \alpha(C^\tau_A; \rho_1,\ldots,\rho_n) =\sum_{i=1}^n p_i \tilde\rho_i ,
$$
$$ 
\tilde\rho_i =\frac{A_{\gamma_i}\rho_i A_{\gamma_i}^*}{\Tr(A_{\gamma_i}^* A_{\gamma_i}\rho_i)} \ \ \ \text{ and } \ \ 
p_i=\Tr(A_{\gamma_i}^* A_{\gamma_i}) ,
$$
where $A_{\gamma_i}=A_{e_{i,1}}\cdots A_{e_{i,m_i}}$, along the oriented path 
$\gamma_i=e_{i,1},\ldots, e_{i,m_i}$ from the $i$--th
leaf to the root. 
$\blacksquare$

\bigskip
\centerline{\bf 5. OPERADS AND ALMOST-SYMPLECTIC QUANTUM CODES} 

\bigskip

In this section we continue our investigation of operadic structures in quantum information, by
revisiting a construction of quantum codes from (not always Moufang) loops obtained from almost-symplectic 
vector spaces over finite fields that we developed in [CoMaMar21]. Here we need to consider a slightly
different definition of almost-symplectic structure over finite fields with respect to the one used in [CoMaMar21]. 
The choice we consider here is better because it allows for an operadic composition, but the
construction of the associated quantum codes is then less well behaved. We show that the space
of almost-symplectoc structures (in the sense we consider here) is an algebra over 
operad modelled on May's little square operad, and the set of data defining the associated quantum
codes is a partial-algebra over the same operad.

\medskip

{\bf 5.1. Rational and binary little square operads.}
The little square operad [May72] provides a characterization of topological spaces that are $2$-fold
loop spaces. Little $n$-cube operads [May72] similarly characterize $n$-fold loop spaces. 

\smallskip

We consider first a sub-operad where we impose an additional condition on the linearly scaled versions
of the unit square in the operad, namely that they have corners located at rational points in the
unit square, namely that the scaling is affected by linear functions with rational coefficients. 

\smallskip

Let $\cI=[0,1]$ be the unit interval, with $\cJ=(0,1)$ its interior. Let $\cI_\Q=\cI\cap \Q$ and $\cJ_\Q=\cJ\cap \Q$.
Thus, $\cI^2$ is the unit square with its interior $\cJ^2$ and $\cI^2_\Q$ and $\cJ^2_\Q$ are the respective sets
of rational points. A rational little square is a function
$$ c: \cI^2_\Q \to \cI^2_\Q, \ \ \  c=(c^1,c^2), \ \ \ c^i(t)=(y_i-x_i)t +x_i \, \ \ \forall t\in \cI_\Q, $$
for some $x_i,y_i\in \cI_\Q$. An $n$-tuple $\langle c_1,\ldots, c_n \rangle$ of rational little squares
has disjoint interiors if $c_i(\cJ^2_\Q)\cap c_j(\cJ^2_\Q)=\emptyset$ for $i\neq j$.

\smallskip

The $n$-th object of the rational little square operad $\cC_2^\Q(n)$ is the space of $n$-tuples
$\langle c_1,\ldots, c_n \rangle$ of rational little squares with disjoint interiors. For $n=0$,
the space $\cC_2^\Q(0)$ consists of a unique function $\emptyset \to \cI^2$.

\smallskip

Let $\sqcup^n \cI^2$ denote the disjoint union $\cI^2\sqcup\cdots\sqcup \cI^2$ of $n$ copies of $\cI^2$.
By identifying $\langle c_1,\ldots, c_n \rangle$ with a function 
$c_1\sqcup\cdots \sqcup c_n: \sqcup^n \cI^2 \to \cI^2$, the set $\cC_2^\Q(n)$
is endowed with the topology induced by the compact-open topology on the set of
continuous maps from $\sqcup^n \cI^2$ to $\cI^2$. 

\smallskip

The operad compositions
$$ \circ_i : \cC^\Q_2(n) \times \cC^\Q_2(m) \to  \cC^\Q_2(n+m-1) $$
of $c=\langle c_1,\ldots, c_n\rangle \in  \cC^\Q_2(n)$ and $c'=\langle c'_1,\ldots, c'_m\rangle \in  \cC^\Q_2(m)$
are determined by the diagrams
$$ \diagram \sqcup^{n+m-1} \cI^2_\Q \ar[rr]^{c\circ_i c'} \ar[dr]_{id_{i-1}\sqcup c' \sqcup id_{n-i}} & & \cI^2_\Q   \\  & \sqcup^n \cI^2_\Q \ar[ur]_c & \enddiagram $$

\smallskip

The unit of the rational little square operad is the identity map $id: \cI^2_\Q \to \cI^2_\Q$ in $\cC_2^\Q(1)$.
The action on $c=\langle c_1,\ldots, c_n \rangle$ of a permutation $\sigma$ in the symmetric group $\Sigma_n$ is given by 
$$ c\sigma := \langle c_{\sigma(1)},\ldots, c_{\sigma(n)} \rangle \, , $$
permuting the labels of the little squares. 

\smallskip
{\bf 5.1.1. Binary little square operad.}
We then consider a subspaces $\cC^{\F_2}_2(n)\subset \cC_2^\Q(n)$ given by those
rational little squares $c=\langle c_1,\ldots, c_n \rangle$ with disjoint interiors, with
the property that the endpoints of each $c_i(\cI^2)$ are rational points of $\cI^2$ that
lie on the square grid of length $2^{-N}$, for some $N\geq 0$.

\smallskip

We refer to the parallel grid of length $2^{-N}$ in the unit square $\cI^2$ as {\it the $N$-grid}.

\smallskip

{\bf 5.1.2. Lemma.}
{\it The subspaces $\cC^{\F_2}_2(n)\subset \cC_2^\Q(n)$, with the induced
composition operations, determine an operad  $\cC^{\F_2}_2$, the
``binary little square operad".}

\smallskip

{\it Proof.} Under the operad composition
operations of $\cC_2^\Q$, the compositions $c\circ_i c'$ of two binary little squares is
still a binary little square, so $\cC^{\F_2}_2$ is a sub-operad of $\cC_2^\Q$.
$\blacksquare$

\smallskip
{\bf 5.1.3. Strict binary little squares.}
Given a binary little square $c\in \cC^{\F_2}_2(n)$, with $c=\langle c_1,\ldots, c_n \rangle$, let $N_c\in \N$ 
be the smallest natural number such that the corners of all the $c_i(\cI^2)$ are at vertices of the $N_c$-grid of
size $2^{-N_c}$. 

\smallskip

{\bf 5.1.4. Definition.} {\it 
A  binary little square $c$ is ``strict" if every row and 
column of the $N_c$-grid has at least one square
that is not contained in the union $\cup_{i=1}^n c_i(\cI^2)$. }

\smallskip

Consider the sub-spaces
$\cC^{\F_2,s}_2(n)\subset \cC^{\F_2}_2(n)$ consisting of binary 
little squares that are strict.

\smallskip

{\bf 5.1.5. Lemma.} {\it The $\cC^{\F_2,s}_2(n)$ with the induced
composition operations, determine a sub-operad of $\cC^{\F_2}_2$.}

\smallskip

{\it Proof.} We need to check that the operad compositions preserve the strict property of
little squares. Given $c\in \cC^{\F_2,s}_2(n)$ and $c'\in \cC^{\F_2,s}_2(m)$, 
the endpoints of the regions $(c\circ_i c')_j(\cI^2)$, $j=1,\ldots, n+m-1$ are on a grid of size
$2^{-N_{c\circ_i c'}}$, with $N_{c\circ_i c'}\geq N_c$. 
Suppose there is a row $R$ (or column) of the $N_{c\circ_i c'}$-grid 
that is completely contained
in the region $\cup_{j=1}^{n+m-1} (c\circ_i c')_j(\cI^2)$. If $R$ does not intersect the region
$c_i(\cI^2)$, then it is contained in the union of the $c_j(\cI^2)$ with $j\neq i$. Since this region
has all sides along the $N_c$-grid, this implies that there must be in fact a row of the $N_c$-grid,
containing $R$, that is contained in $\cup_{j\neq i} c_j(\cI^2)$, but this is not possible because
$c$ is strict. Thus, $R$ must intersect the region $c_i(\cI^2)$. This means that, within the region
$c_i(\cI^2)$ a row of the $N_{c\circ_i c'}$-grid is completely contained in the union of the
linearly scaled images of the $c'_j(\cI^2)$, but this in turn implies that in $\cI^2$ a row of
the $N_{c'}$-grid must be contained in the union of the $c'_j(\cI^2)$, which cannot happen
because $c'$ is strict.
$\blacksquare$

\medskip

{\bf 5.2. Binary little square operads and almost symplectic spaces.}
We consider here a class of (not always Moufang) loops that we previously introduced and investigated in [CoMaMar21]. These are
obtained from almost-symplectic structures on vector spaces over a finite field of characteristic $2$. We recall the basic
setting from [CoMaMar21]. However, the notion of almost-symplectic form we consider here is somewhat different from
[CoMaMar21]: this will allow for better properties with respect to operadic composition, but will in turn have worse properties
with respect to representations of the resulting loops, hence the construction of quantum codes considered in [CoMaMar21]
will not directly extend to this setting.

\smallskip

Let $\F=\F_{2^r}$ be a finite field of characteristic two and let $V$ be a finite dimensional vector space over $\F$. 
Let $K$ be an unramified extension of $\Q_2$ with residue field $\cO_K/\fm_K=\F$. Consider the ring $R=\cO_K/\fm_K^2$
and a free $R$-module $\tilde V$ with $V=\tilde V/\fm_K$. Consider functions $\tilde\omega: \tilde V \times \tilde V\to R$
with $\omega =2 \tilde\omega$ the induced function $\omega: V\times V \to R$. Note that these functions are not
linear as they do not satisfy the Hochschild cocycle condition of symplectic forms, namely we require that
$$ d\omega(u,v,w)=\omega(v,w) -\omega(u+v,w) +\omega(u,v+w) -\omega(u,v) \neq 0\, . $$

\smallskip

A pair $(V,\omega)$ as above is an almost-symplectic space if $\omega$ is non-degenerate, in the sense that
for any $u \in V$ there is some $v\in V$ with $(u,v)\neq (0,0)$, such that $\omega(u,v)\neq 0$. 
(Note that with this definition the almost-symplectic structure $\omega$ does
{\it not} satisfy $\omega(0,v)=\omega(u,0)=0$, for all $u,v\in V$.)
The almost-symplectic structure $\omega$ has a polarization $\beta: V \times V \to R$
satisfying $\beta(u,v)-\beta(v,u)=\omega(u,v)$. 
The loop $\cL(V,\beta)$ associated to the data $(V,\beta)$ is the 
extension
$$ 0 \to R \to \cL(V,\beta) \to V \to 0 $$
with non-associative multiplication
$$ (x,u)\star (y,v)=(x+y+\beta(u,v), u+v) \, . $$
The non-associativity is a consequence of the fact that $(V,\omega)$ is almost-symplectic and not symplectic,
hence the Hochshild coboundary of $\beta$ is nonzero, 
$$ d\beta(u,v,w)=\beta(v,w) -\beta(u+v,w) +\beta(u,v+w) -\beta(u,v) = \gamma(u,v,w) \not\equiv 0 \, . $$

\smallskip

The Moufang condition for the loop $\cL(V,\beta)$ is equivalent to an identity for $\gamma =d\beta$, see [CoMaMar21]. 
We do not necessarily require here the condition that the loop $\cL(V,\beta)$ is Moufang. 

\smallskip

We focus in particular on the special case where $\F=\F_2$ and $R=\Z/4\Z$. 
We can identify the map $\omega=2\tilde\omega: V \times V \to R$, where $V=\F_2^N$, as a subdivision of the 
square $\cI^2=[0,1]\times [0,1]$ into the $N$-grid of $2^N\times 2^N$ sub-squares of side $2^{-N}$, which we can label
by the pairs $(u,v)\in V\times V$. Each square is colored white or black according to whether $\omega(u,v)=0$ or not. 
The property that $\omega$ is almost-symplectic (that is, non-degenerate) is equivalent to the fact that in every
row and column of the subdivided square there is at least one black sub-square. 

\smallskip

Let  $\cV^{\F_2}$ be the space of almost symplectic finite dimensional vector spaces over $\F_2$,
topologized as a subset of the space of maps $\cup_N Maps(\F_2^N\times \F_2^N, R)$, where
$\cV^{\F_2}$ consists of those maps that satisfy the non-degenerate condition above, with 
range $2R\subset R$, since $\omega =2 \tilde\omega$. 

\smallskip

{\bf 5.2.1. Theorem.} {\it The space $\cV^{\F_2}$ of almost symplectic finite dimensional vector spaces over $\F_2$
is an algebra over the operad $\cC^{\F_2,s}_2$.}

\smallskip

{\it Proof.}
To realize $\cV^{\F_2}$ as an algebra over $\cC^{\F_2,s}_2$ we need operations
$$ \cC^{\F_2,s}_2(n) \times\underbrace{\cV^{\F_2}\times \cdots \times\cV^{\F_2}}_{n-\text{times}} \to \cV^{\F_2}, $$
that assign to an $n$-tuple $(V_1,\omega_1), \ldots, (V_n, \omega_n)$ and a
strict binary little square $c=\langle c_1,\ldots, c_n \rangle$ a new
almost-symplectic space $(V,\omega)=\gamma(c; (V_1,\omega_1), \ldots, (V_n, \omega_n))$,
compatibly with the composition operations in the operad $\cC^{\F_2,s}_2$. 
We construct $(V,\omega)$ in the following way. For $V_i=\F_2^{N_i}$, we consider
the regions $\omega_i^{-1}(0)$ with sides on the $N_i$-grid in $\cI^2$. We take each
of these copies of $\cI^2$ subdivided in the the subsquares of the $N_i$-grid, with
those in $\omega_i^{-1}(0)$ colored white and the others colored black, and we scale
it linearly so as to fit, respectively, into the $c_i(\cI^2)$ regions of $\cI^2$ determined 
by the binary little square $c=\langle c_1,\ldots, c_n \rangle$. We color black the outside of
$\cup_i c_i(\cI^2)$ in $\cI^2$. Let $N\in \N$ be the smallest integer such that all the
resulting contours in $\cI^2$ separating the black and the white colored areas are
on the $N$-grid. We then set $V=\F_2^N$ and we assign the values of $\omega(v,w)$
according to the color of the corresponding square. The fact that the binary little square
is strict implies that $\omega$ is non-degenerate. 
$\blacksquare$

\medskip
{\bf 5.3. Colored $p$-ary little squares.} Let now $q=p^r$ be some prime power with $p>2$. We construct an operad
that generalizes the strict binary little squares operad considered in the case of characteristic $2$.

\smallskip

For simplicity, as in the case of characteristic $2$, we restrict to the case of $\F_p$. We then define the $p$-ary $N$-grid
in the unit square $\cI^2$ to the the parallel grid of size $p^{-N}$, with $p^N\times p^N$ sub-squares. We just refer to it
as the $N$-grid when the choice of $p$ is understood.

\smallskip

{\bf 5.3.1. Definition.} {\it 
A colored $p$-ary little square is a decomposition of the unit square $\cI^2$ into $p$ regions $\cR_0,\ldots, \cR_{p-1}$
with disjoint interiors $\cJ(\cR_i)\cap \cJ(\cR_j)=\emptyset$, for $i\neq j$, with $\cJ(\cR)$ denoting the interior of a region $\cR$. 
With the property that, for each $i=0,\ldots, p-1$, and for some integers $n_i, N_i\in \N$, there is a rational little square $c^{(i)}=\langle c^{(i)}_1,\ldots, c^{(i)}_{n_i} \rangle$
in $\cC_2^\Q(n_i)$, with endpoints on a $p$-ary $N_i$-grid, such that $\cR_i=\cup_{j=1}^{n_i} c^{(i)}_j(\cI^2)$. 
The colored $p$-ary little square is strict if,  for $i=0$, the little square $c^{(0)}$ is strict, namely no row or column of the $p$-ary 
$N_0$-grid is completely contained
in $\cup_{j=1}^{n_0} c^{(0)}_j(\cI^2)$. We denote by $\cC^{\F_p}_2(n_0,\ldots,n_{p-1})$ the set of colored $p$-ary little squares as above
and we take
$$ \cC^{\F_p}_2(n):= \bigcup_{n_1,\ldots, n_{p-1}} \cC^{\F_p}_2(n, n_1, \ldots,n_{p-1})\, . $$
The set of strict colored $p$-ary little squares $\cC^{\F_p,s}_2(n)$ is similarly defined. }

\smallskip

We denote the (strict) colored $p$-ary little squares in $\cC^{\F_p}_2(n)$ (or $\cC^{\F_p,s}_2(n)$, respectively)  with the notation 
$$ \cI^2(\cR_0,c^{(0)};\cR):=\left\{ (\cR_0,c^{(0)}=\langle c^{(0)}_1,\ldots, c^{(0)}_n \rangle), (\cR_j,c^{(j)})_{j=1,\ldots, p-1}\right\}  \, , $$
with $\cR:=\{ (\cR_j,c^{(j)}) \}_{j=1,\ldots, p-1}$.
The composition operations 
$$ \circ_i: \cC^{\F_p}_2(n) \times \cC^{\F_p}_2(m) \to \cC^{\F_p}_2(n+m -1) $$
are determined by taking
$$ \cI^2(\cR_0,c^{(0)};\cR)\circ_i \cI^2(\cR'_0,c^{\prime (0)};\cR') $$
to be given by the decomposition of the unit square $\cI^2$ into the regions 
$$ \cI^2 = \cB \cup \bigcup_{j\neq 0} \cR_j \cup \bigcup_{\ell\neq i} c^{(0)}_\ell (\cI^2) \, , $$
where the region $\cB$ is obtained by linearly scaling the square $\cI^2=\cup_j \cR^\prime_j$ subdivided into the regions of the
second colored $p$-ary little square and placing it in place of the region $c^{(0)}_i(\cI^2)$. 

\smallskip

\smallskip

{\bf 5.3.2. Lemma.} {\it
With the composition operations above the $\cC^{\F_p}_2(n)$ (respectively,  $\cC^{\F_p,s}_2(n)$) form an operad. }

\smallskip

{\it Proof.}
The new regions $\cR''_\ell$ with $\ell=1,\ldots,p-1$ of the composed $\cI^2(\cR_0,c^{(0)};\cR)\circ_i \cI^2(\cR'_0,c^{\prime (0)};\cR')$
are given by
$$ \cR''_\ell =\cR_\ell \cup c^{\prime (0)}_i(\cR'_\ell)\, $$
while the $\cR''_0$ region of the composed $\cI^2(\cR_0,c^{(0)};\cR)\circ_i \cI^2(\cR'_0,c^{\prime (0)};\cR')$ is given by
$$ \cR''_0 = \bigcup_{j\neq i} c^{(0)}_j(\cI^2) \cup c^{\prime (0)}_i(\cR'_0)\, . $$
Thus we see that these composition operations
are still the same composition operations of the operad $\cC_2^\Q$, acting on the $c^{(0)}$ little squares while maintaining all the 
rest of the data unaffected. The strict condition is preserved under composition by the same argument as in Lemma~5.1.5. 
$\blacksquare$

\medskip
{\bf 5.3.3.  Operads and almost-symplectic structures over $\F_p$.}
An almost-symplectic vector space over $\F_p$ is a pair $(V,\omega)$ of a finite dimensional vector
space over $\F_p$ and a function $\omega: V\times V \to \F_p$ which is non-degenerate, in the sense that
for all $u\in V$ there is some $v$ with $(u,v)\neq (0,0)$ such that $\omega(u,v)\neq 0$. and with nontrivial
Hochschild coboundary 
$$ d\omega(u,v,w)=\omega(v,w) -\omega(u+v,w) +\omega(u,v+w) -\omega(u,v) \neq 0\, . $$
In this case again the non-vanishing of the Hochschild coboundary  implies that $\omega$ cannot be
a linear map. 

\smallskip

Note that here also the notion of almost symplectic structure we are considering is
different from [CoMaMar21], hence the construction of representations and quantum codes described 
there does not apply directly to this case.

\smallskip

In this case, the (not necessarily Moufang) loop $\cL(V,\omega)$ associated to the almost symplectic structure
$(V,\omega)$ is obtained as the extension
$$ 0 \to \F_p \to \cL(V,\omega) \to V \to 0 $$
with the non-associative multiplication given by
$$ (x,u)\star (y,v)=(x+y+\frac{1}{2}\omega(u,v), u+v)\, . $$
As shown in [CoMaMar21], the Moufang condition for $\cL(V,\omega)$ is expressed as an identity satisfied by the function $d\omega$.

\smallskip

Let $\cV^{\F_p}$ denote the space of almost symplectic structures $(V,\omega)$ over $\F_p$.

\smallskip

{\bf 5.3.4. Theorem.} {\it
The space $\cV^{\F_p}$ is an algebra over the operad $\cC^{\F_p,s}_2$ of strict colored $p$-ary little squares. 
}

\smallskip

{\it Proof.} The argument is similar to Theorem~5.2.1, except for the coloring of the regions. Given an $n$-tuple
of almost symplectic spaces $(V_i,\omega_i)$ and a strict colored $p$-ary little square $\cI^2(\cR_0,c^{(0)};\cR)\in \cC^{\F_p,s}_2(n)$,
we form a new $$(V,\omega)=\gamma(\cI^2(\cR_0,c^{(0)};\cR); (V_1,\omega_1),\ldots,(V_n,\omega_n)) $$ by 
associating to each $(V_i,\omega_i)$, with $V_i=\F_p^{N_i}$, a $p$-ary $N_i$-grid in $\cI^2$ subdivided into 
regions $\cR_\ell=\omega_i^{-1}(\ell)$ for $\ell\in \F_p$. This determines a colored $p$-ary little square $\cI^2(\cR_{i,0},c^{i,(0)};\cR_i)$ associated
to each $(V_i,\omega_i)$, which is strict because $\omega_i$ is non-degenerate. We then compose these
according to the composition 
$$ \gamma(\cI^2(\cR_0,c^{(0)};\cR);\cI^2(\cR_{1,0},c^{1,(0)};\cR_1),\ldots, \cI^2(\cR_{n,0},c^{n,(0)};\cR_n)) $$
$$ \gamma: \cC^{\F_p,s}_2(n)\times \cC^{\F_p,s}_2(k_1)\times \cdots \times \cC^{\F_p,s}_2(k_n)) \to \cC^{\F_p,s}_2(k_1+\cdots+k_n) $$
obtained from repeated application of the compositions $\circ_i$ in the operad $\cC^{\F_p,s}_2$. This results in a new
colored $p$-ary little square with regions $\cR_\ell''$, which is also strict. This in turn defines the resulting almost symplectic space $(V,\omega)$ which has
$V=\F_p^N$, with $N$ the smallest natural number such that the subdivision into regions $\cR''_\ell$ of the resulting colored $p$-ary little square
is along the $p$-ary $N$-grid, and $\omega(u,v)=\ell$ for $(u,v)\in \cR''_\ell$. This $\omega$ is non-degenerate because the little square is strict.
$\blacksquare$

\medskip
{\bf 5.4. Operad partial-action on quantum codes.} 
A vector space $V=\F_q^N$, determines a corresponding complex vector space $(\C^q)^{\otimes N}$,
representing a system of $N$ $q$-ary qbits, endowed with the canonical basis given by the $| v \rangle$, labelled by the vectors $v\in V$. 

\smallskip

We consider here the case of $\F_p$ with $p>2$. The argument can be adapted to
the case of characteristic $2$ along the lines discussed in [CoMaMar21]. Let $\cL=\cL(V,\omega)$ be the loop obtained
from an almost-symplectic $(V,\omega)$ over $\F_p$ as recalled above. It acts on $\cH=\C[\cL(V,\omega)]$ by left and
right multiplication $((x,u)\star f)(y,v)=f((x,u)\star (y,v))$ and $(f\star (x,u))(y,v)=f((y,v)\star(x,u))$. 
This gives rise to a loop representation, which is the product $\cL\times \cH$ endowed with the
(non-associative) multiplication
$$ ((x,u),f)\star ((y,v),f')=((x,u)\star(y,v), (x,u)\star f' + f \star (y,v))\, . $$

\smallskip

Consider the following sets
Let $\cS^\omega_1=\omega^{-1}(0)\cap \omega^{-1}(0)^\tau=\{ (u,v)\in V\times V \,|\, \omega(u,v)=\omega(v,u)=0 \}$.
Given a subset $\Omega\subset V\times V$, we write $$\Omega_{ij}=\pi_{ij}^{-1}\Omega \subset V^L=\underbrace{V\times\cdots\times V}_{L-\text{times}}\, , $$
with $\pi_{ij}$ the projection onto the product of the $i$-th and $j$-th components. 
For $L\geq 2$ let $\cS_L \subset V^L$, 
$$ \cS^\omega_L = \bigcap_{1\leq i<j\leq L} (\cS^\omega_1)_{ij}\, . $$
Let $\cS=\cup_{(V,\omega)\in \cV^{\F_p}} \cup_{L\geq 1} \cS^\omega_L$. 

\smallskip

Given a character $\chi: \F_p\to \C^*$, let $\cH_\chi \subset \cH$ be the subspace of functions
$f(x,u)$ that transform by $((x',0)\star f)(x,u)=\chi(x') f(x,u)$. An element $\underline{u}=(u_1,\ldots, u_L)\in \cS_L$
determines a set $\{ \chi(x_1) E_{u_1},\ldots, \chi(x_L) E_{u_L} \}_{x_i\in \F_p}$ of mutually commuting
operators on $\cH_\chi$. This determines a quantum code $\cQ^\lambda_{\chi,\underline{u}}\subset \cH$ given by
a common eigenspace of these operators with eigenvalue $\lambda$. We refer to the quantum codes obtained in this
way as ``almost-symplectic quantum codes".

\smallskip

The operad action on almost-symplectic structures described in Theorem~5.3.4 induces a partial-action of the same operad
on the data that determine these quantum codes. The notion of partial-action of an operad
(partial-algebra over an operad) was introduced in [KrMay95].

\smallskip

\smallskip

{\bf 5.4.1. Proposition.} {\it The space $\cS$ is a partial-algebra over the operad $\cC^{\F_p,s}_2$ of strict colored $p$-ary little squares. }

\smallskip

{\it Proof.} In $\cS^n=\underbrace{\cS\times\cdots\times \cS}_{n-\text{times}}$ consider the subspace 
$$ \cS^n_0=\{ (\underline{u}^{(1)},\ldots,\underline{u}^{(n)})\in \cS^{\omega_1}_{L_1}\times \cdots \times \cS^{\omega_n}_{L_n} \,|\,  
(\tilde{\underline{u}}^{(1)},\ldots,\tilde{\underline{u}}^{(n)}) \in \cS^\omega_{L_1+\cdots + L_n}   \}\, $$
for $(V,\omega)=\gamma(\cI(\cR_0,c^{(0)};\cR); (V_1,\omega_1),\ldots, (V_n,\omega_n))$ the operad action of Theorem~5.3.4,
and with $(\tilde{\underline{u}}_i)$ the vector in $V$ obtained from the vectors $\underline{u}_i \in V_i$ under this composition.
Then the operad action on $\cV^{\F_p}$ induces a partial-action
$\gamma: \cC^{\F_p,s}_2(n)\times \cS^n_0 \to \cS$.
$\blacksquare$ 

\bigskip
\bigskip

{\bf Acknowledgment.} The first author is supported by the Minerva Fast track grant from the Max
Planck Institute for Mathematics in the Sciences. The third author is supported by NSF grants 
DMS-1707882 and DMS-2104330.

\bigskip
\bigskip

\centerline{\bf References}

\medskip

[BFL11] J.C.~Baez, T.~Fritz, T.~Leinster, {\it A characterization of entropy in terms of information loss}, 
Entropy 13 (2011), no. 11, 1945--1957. 

\smallskip

[BoMa08]  D.~Borisov, Yu. Manin. {\it Generalized operads and their inner cohomomorhisms.}
 In: Geometry and Dynamics of Groups
and Spaces (In memory of Aleksander Reznikov). Ed. by M. Kapranov et al.
Progress in Math., vol. 265 (2008), 
Birkh\"auser, Boston, pp. 247--308.
arXiv:math.CT/0609748.

\smallskip

[ChCorGi19]  C.  Chenavier,  C.  Cordero,  S.  Giraudo.  {\it Quotients of the magmatic
operad: lattice structures and convergent rewrite systems.} arXiv:1809.05083v2, 30 pp.

\smallskip

[CoMa21] N. Combe, Yu. Manin. {\it Symmetries of genus zero modular operad.} In ``Integrability, Quantization, 
and Geometry: II. Quantum Theories and Algebraic Geometry", Proceedings of Symposia in Pure Mathematics
Vol.103, pp.101--110, American Mathematical Society, 2021.
arXiv:math.AG/1907.10317.

\smallskip

[CoMaMar21]  N. C. Combe, M. Marcolli, Yu. Manin. {\it Moufang patterns and geometry of information.}
To be published in Collection, dedicated to Don Zagier, in Pure and Applied Math. Quarterly.
arXiv:2107.07486. 42 pp.

\smallskip

[DeMi82] P. Deligne,  J.  Milne.  {\it Tannakian categories.} In ``Hodge cycles, motives and Shimura 
varieties, Springer Lecture Notes in Math. 900 (1982), pp. 101--228.

\smallskip

[KrMay95] I.~K\v{r}\'i\v{z}, J.P.~May, 
{\it Operads, algebras, modules and motives}. Ast\'erisque No. 233 (1995), iv+145pp.

\smallskip

[Ma87] Yu.  Manin.  {\it Some remarks on Koszul algebras and quantum groups}.  Ann. Inst.
Fourier, XXXVII (1987) no.  4, pp. 191--205.

\smallskip

[Ma88]  Yu.  Manin.  {\it Quantum groups and non--commutative geometry}. Publ. de CRM,
Universit\'e  de Montr\'eal (1988), 91 pp.

\smallskip

[Ma91]  Yu. Manin.  {\it Topics in noncommutative geometry}. Princeton University Press (1991), 163 pp.

\smallskip

[Ma92]  Yu.  Manin.  {\it Notes on quantum groups and quantum de Rham complexes.}
Teoreticheskaya  i Matematicheskaya Fizika 92:3 (1992), pp. 425--450.  Reprinted
in Selected Papers of Yu.  I.  Manin, World Scientific, Singapore 1996, pp. 529--554.

\smallskip

[MarThor14] M.~Marcolli, R.~Thorngren, {\it Thermodynamic semirings}, 
J. Noncommut. Geom. 8 (2014), no. 2, 337--392

\smallskip

[Markl08] M.~Markl, {\it Operads and PROPs},  Handbook of algebra. Vol. 5, 87--140,
Elsevier/North-Holland, 2008.

\smallskip

[May72] J.P.~May, {\it The geometry of iterated loop spaces}, Lecture Notes in Mathematics, Vol.~271, Springer, 1972.

\smallskip

[Sm16] J. D. H. Smith. {\it Quantum quasigroups and loops.}  Journ. of Algebra 456 (2016),
pp. 46--75.

\smallskip

[Val04] V.  Vallette. {\it A Koszul duality for PROPs.} arXiv:math0411542v3, 78 pp.

\smallskip

[Val08] B. Vallette.  {\it Manin products, Koszul duality, Loday algebras and Deligne conjecture.}
J. Reine Angew. Math. vol. 620 (2008), pp. 105 -- 154. arXiv:math/0609002.

\smallskip

[Vor01] A.~Voronov, {\it The $A_\infty$ operad and $A_\infty$-algebras}, 
lecture notes, 2001. \newline {\tt https://www-users.cse.umn.edu/~voronov/8390/lec9.pdf}

\enddocument